\documentclass[11pt,reqno]{amsart}
\usepackage[cp1251]{inputenc}
\usepackage[english]{babel}
\usepackage{amsmath}
\usepackage{amssymb}
\usepackage{amsfonts}
\usepackage{xcolor}
\usepackage[colorlinks]{hyperref}
 \usepackage{amsmath,systeme}
\usepackage{nicefrac}

\usepackage{titlesec}
\titlespacing{\section}{0pt}{2.5ex}{1.5ex}
\titlespacing{\subsection}{0pt}{1.5ex}{1ex}
\titlespacing{\subsubsection}{0pt}{1.5ex}{1ex}
\titleformat{\section}{\large\bfseries\centering}{\thesection}{1em}{}
\titleformat{\subsection}[runin]{\bfseries}{\thesubsection.}{0.5em}{}[.\mbox{\ }]
\titleformat{\subsubsection}[runin]{\bfseries}{\thesubsubsection.}{0.4em}{}[.\mbox{\ }]

\newtheorem{remark}{Remark}

\newtheorem{theorem}{Theorem}
\newtheorem{definition}{Definition}

\newtheorem{example}{Example}
\newtheorem{problem}{Problem}

\definecolor{sangria}{rgb}{0.57, 0.0, 0.04}

\hypersetup{linkcolor=red,filecolor=black,urlcolor=blue, citecolor=blue}

\begin{document}
\renewcommand{\refname}{References}
\renewcommand{\proofname}{Proof.}
\renewcommand{\figurename}{Fig.}

\thispagestyle{empty}

\title[Special Pythagorean Triangle]{\large A Pythagorean triangle in which the hypotenuse and the sum of the arms are squares}

\author[Djamel HIMANE]{{Djamel Himane}   }

\address{Djamel Himane   
\newline\hphantom{iii} LA3C Laboratory,  
\newline\hphantom{iii} Faculty of Mathematics, USTHB, 
\newline\hphantom{iii} Po. Box 32, El Alia, 16111, 
\newline\hphantom{iii} Bab Ezzouar, Algiers, Algeria.}%
\email{\textcolor{blue}{dhimane@usthb.dz}}%

 \maketitle {\small

\vspace{-12pt}
 
\bigskip
\bigskip

\begin{quote}
\noindent{\bf Abstract:} In this paper, show that the Diophantine equation $ x^2+(x+1)^2=w^4 $ has only two solutions $ (0,1) $ and $  (119,13)$ in non-negative integers $ x $ and $ w $. This equation concerned a classic problem posed by Pierre de Fermat, wonders about finding a Pythagorean triangle in which the hypotenuse and the sum of the arms are square. We review the method of finding the smallest  solution presented by Fermat, and the relationship between the primitive Pythagorean triples and the Pell's equation, Finally, we present an algorithm for finding primitive solutions, which actually enabled us to find other solutions.

 \end{quote}
}

\bigskip

\section{Introduction}

A right triangle is a triangle whose sides $(x, y, z)$ satisfy the Diophantine equation
\begin{equation} \label{eq.01} 
 x^2 + y^2 = z^2 
\end{equation}
A primitive triangle is an integral triangle such that the greatest common divisor
of the lengths of its sides is 1. A right triangle is a Pythagorean triangle if the right triangle has integral sides. A primitive Pythagorean triangle has sides
 \begin{equation} \label{eq.02}
 (x, y, z)= (m^2-n^2, 2mn ,  m^2+n^2 )
 \end{equation}
 where
 \begin{equation*}
  m > n > 0, \hspace{1cm} \gcd(m, n) = 1, \hspace{1cm}  m + n \equiv 1 \pmod 2 .
 \end{equation*}
 
In 1643, Fermat stated the following problem: (See page 67 in \cite{1})
\\

\textit{Find a Pythagorean triangle in which the hypotenuse and the sum of the arms are squares.}
\\

In a letter to Mersenne, Fermat affirmed that the smallest of these triangles is the triangle
  \begin{equation} \label{eq.03} 
  (4565486027761, 1061652293520, 4687298610289). 
\end{equation}
   That is,  
   \begin{equation*}  
  (4565486027761)^2 + (1061652293520)^2 = (4687298610289)^2 = (2165017)^4
\end{equation*}
and
\begin{equation*}  
   4565486027761  +  1061652293520  = (2372159)^2 
\end{equation*}
\\

The method of Fermat is presented in problems $ n^{\circ} $ 22, 23, pp. 339-340, and $ n^{\circ} $ 45, pp. 453-354 (see \cite{2}).  The solution presented above was reformulated by W. Sierpinski \cite{1} as follows:
\\

First of all we show how one can arrive by the shortest method at the solution given by Fermat.

The problem posed by Fermat is obviously equivalent to the problem of obtaining the solutions $x, y, e, f  $ of the simultaneous equations
 \begin{equation} \label{eq.04}
 x^2 + y^2 = e^4, \hspace{2cm} x + y = f^2
 \end{equation}
in natural numbers. It is sufficient to obtain the solutions of the equations \eqref{eq.04} in positive rational numbers $ x, y, e, f,$ because reducing them to a common denominator $ h $ and multiplying the first of the equations \eqref{eq.04} by $ h^4 $ and the second by $ h^2 $, we obtain the solutions of the equations \eqref{eq.04} in natural numbers $ h^2 x, h^2 y, h e, h f $.
We suppose for a rational number $t$,
 \begin{equation} \label{eq.05}
 x = (t + 5)^2 - 12^2, \hspace{2cm} y = 24(t + 5).
 \end{equation}
We have the identity
 \begin{equation*}
\left[(t + 5)^2 - 12^2 \right]^2 + \left[24(t + 5)\right]^2 = \left[(t + 5)^2 + 12^2\right]^2
  \end{equation*}
In order that the numbers \eqref{eq.05} satisfy equation \eqref{eq.04} it is necessary and sufficient that
 \begin{equation*}
(t + 5)^2 + 12^2 = e^2, \hspace{2cm} (t + 5)^2 - 12^2 + 24(t + 5) = f^2
  \end{equation*}
That is,
\begin{equation} \label{eq.06}
t^2 + 10t + 169 = e^2
 \end{equation}
and
\begin{equation} \label{eq.07}
t^2 + 34t + 1 = f^2
 \end{equation}
The equality \eqref{eq.07} gives
\begin{equation*} 
 169(t^2 + 34t + 1) = (13f)^2,
  \end{equation*}
which because of \eqref{eq.06} gives
\begin{equation} \label{eq.08}
 (13f)^2 - e^2 = 168t^2 + 57364 t= 14t \left( 12t + \dfrac{2868}{7} \right) 
  \end{equation}
The equality \eqref{eq.08} will hold if
\begin{equation*} 
13 f-e = 14t, \hspace{2cm} 13f + e = 12t + \dfrac{2868}{7} 
  \end{equation*}
hence
 \begin{equation} \label{eq.09}
  e = -t + \dfrac{1434}{7}  \hspace{1cm} \text{and} \hspace{1cm} f = t + \dfrac{1434}{7.13}.  
  \end{equation}
The first of the formula \eqref{eq.09} gives
\begin{equation} \label{eq.10}
  e^2 = t^2 - \dfrac{2868}{7} t + \dfrac{1434^2}{7^2}. 
  \end{equation}
In order, therefore, that the number $ e^2 $ satisfy equation \eqref{eq.06} it is necessary and sufficient that $t$ satisfy the equation
\begin{equation} \label{eq.11}
  10t + 169 = - \dfrac{2868}{7} t  + \dfrac{1434^2}{7^2},
  \end{equation}
that is,
\begin{equation} \label{eq.12}
t= \dfrac{1434^2 - (13.7)^2}{7^2 \left( 10+ \dfrac{2868}{7}\right)} = \dfrac{1343.1525}{7.2938} = \dfrac{2048075}{20566}
  \end{equation}
On the other hand, if we take the value of $t$ from formula \eqref{eq.12} and the values $e$ and $f$ from \eqref{eq.09}, then formula  \eqref{eq.10} and \eqref{eq.11} give formula \eqref{eq.06} and because of the second of the formula  \eqref{eq.09} we obtain formula \eqref{eq.08} which from \eqref{eq.06} gives formula \eqref{eq.06}. But as we know, formula  \eqref{eq.05}, \eqref{eq.06} and \eqref{eq.07} give formula \eqref{eq.04}.

If then we determine the number $t$ from formula \eqref{eq.12} and the numbers $ x, y, e, f $ from formula  \eqref{eq.05} and \eqref{eq.09}, we obtain the solution of the set of equations \eqref{eq.04} in positive rational numbers in which the denominator of each of the numbers $ y, e, f $ will be, as is easy to establish, $ h = 20566 $ and the denominator of the number $ x $ will be $ h^2 $. The numbers 
$ h^2x, h^2y, $ and $ h^2 e^2 $ are then precisely the sides of the pythagorean triangle \eqref{eq.03} given by Fermat.

Thus there exists a pythagorean triangle in which the hypotenuse and the sum of the arms are squares of natural numbers. If $ (x, y, z) $ is such a triangle and $ k $ is an arbitrary natural number, then obviously the triangle $ (k^2 x, k^2 y, k^2 z) $ will also be such a triangle.
\begin{remark}
 The trivial solutions of system \eqref{eq.04} are the two solutions 
\begin{equation*}
 (x,y,e,f)=\lbrace(0,0,0,0),(1,0,1,1)\rbrace.
\end{equation*}
\end{remark}
 All primitive solutions discovered so far can be found in sequences 
\href{https://oeis.org/A166930}{A166930}, \href{https://oeis.org/A166929}{A166929} and \href{https://oeis.org/A167437}{A167437}.
\section{Preliminaries}

We recall \cite{4} some properties of solutions to the negative Pell's equation
 \begin{equation} \label{eq.13}
 u^2-2v^2= - 1 
 \end{equation}
 and the positive Pell's equation
 \begin{equation} \label{eq.14}
 u^2-2v^2= + 1 
 \end{equation}
Suppose that $ (u,v)=(u_{t},v_{t}) $  is a solution to \eqref{eq.13} or \eqref{eq.14}  . Define the
integer pair $(u_{2t},v_{2t}) $ by the equation
\begin{equation*}
    u_{2t}+ v_{2t} \sqrt{2}  = \left( u_{t} + v_{t} \sqrt{2}\right)^2  
\end{equation*}
We can easily check that
\begin{equation} \label{eq.15} 
  u_{2t}= u_{t}^2+2 v_{t}^2, \hspace{2cm} v_{2t} =2u_{t}v_{t} 
\end{equation}
More generally, suppose that $ (u_{k}, v_{k}) $ is defined by
\begin{equation} \label{eq.16} 
    u_{k}+ v_{k} \sqrt{2}  = \left( 1 +  \sqrt{2}\right)^k  \text{for }  k\geq 1.
\end{equation}

thus, for $ t\geq 1 $, then $ (u_{2t-1}, v_{2t-1}) $ is a solution of \eqref{eq.13} and $ (u_{2t}, v_{2t}) $ is a solution of   \eqref{eq.14}. From \eqref{eq.16} , note that
\begin{equation*}
  u_{k}+v_{k}\sqrt{2}= (u_{k-1} +  v_{k-1} \sqrt{2})(1+\sqrt{2}) 
\end{equation*}
  and deduce that
\begin{equation} \label{eq.17}
    u_{k}=u_{k-1}+2 v_{k-1}, \hspace{2cm}  v_{k}  =   u_{k-1} + v_{k-1}  
\end{equation}
 The equality \eqref{eq.17} gives
\begin{equation} \label{eq.18}
    u_{k-1}= 2 v_{k} - u_{k}, \hspace{2cm}  v_{k-1}  =   u_{k} - v_{k}  
\end{equation}
then formula \eqref{eq.15} and \eqref{eq.18} give
\begin{equation} \label{eq.19}
    u_{2t-1}=  u_{t}^2 - 2(u_{t}- v_{t})^2, \hspace{2cm}  v_{2t-1}  =   (u_{t} - v_{t})^2  + v_{t}^2
\end{equation}
\\

R. Guy, in his book of unsolved problems \cite{5}, shows that the Diophantine equation 
\begin{equation*}
x^2 - 2y^4= - 1
\end{equation*}
  has, in the set of positive integers, the only solutions $ (1,1) $ and $ (239,13) $, This equation has been solved in several a method from Ljunggren's, Mordell's and F. Smarandache. (see \cite{3})
\\

These two solutions \cite{3} can be derived through formula \eqref{eq.16} as follows:   
\begin{equation} \label{eq.20} 
    u_{1}+ v_{1} \sqrt{2}  =  1 +  \sqrt{2}  \hspace{.8cm} \text{and} \hspace{.8cm} u_{7}+ v_{7} \sqrt{2}  = 239 + 169 \sqrt{2}
\end{equation}

\section{Main results}

\begin{theorem}
The Diophantine equation
\begin{equation} \label{eq.21} 
  x^2+(x+1)^2 = w^4  
\end{equation}
  has only two solutions $(0, 1) $ and $(119, 13)$  in positive integers $ x, w $.
\end{theorem}

\begin{proof} 
we assume that the Diophantine equation \eqref{eq.21}  has solution, then from \eqref{eq.01}  and formula \eqref{eq.02} there exist positive integers $ m, n, x, z,$ such that 
\begin{equation*}
x=m^2-n^2,  \hspace{1cm}   x+1=2mn, \hspace{1cm}  z=w^2=m^2+n^2
\end{equation*}
or 
\begin{equation*}
x+1=m^2-n^2,  \hspace{1cm}   x=2mn, \hspace{1cm}  z=w^2=m^2+n^2,
\end{equation*}
where 
\begin{equation*}
  m > n > 0, \hspace{.3cm} \gcd(m, n) = 1, \hspace{.3cm}  m + n \equiv 1 \pmod 2  
 \end{equation*}
thus, 
\begin{equation} \label{eq.22}
 m^2-n^2-2mn = \pm1,  
  \end{equation}
 as 
 \begin{equation*}
 \gcd(x, x+1)=1
  \end{equation*}
The equality  \eqref{eq.22} gives
\begin{equation} \label{eq.23}
 (m+n)^2-2m^2 = \mp 1,  
  \end{equation}
by formula \eqref{eq.16}, we obtain 
\begin{equation*}
  u_{s} + v_{s}\sqrt{2} = (m+n) + m \sqrt{2}
  \end{equation*}
 for some positive integers $ s $.
The second of the formula \eqref{eq.19} gives
\begin{equation} \label{eq.25} 
      v_{2s-1}  =   (u_{s} - v_{s})^2  + v_{s}^2 = m^2 + n^2 = z
\end{equation}
The equality  \eqref{eq.20} gives $ s=1 $ or $ s=4 $, for $ s=1 $, We get the trivial solution
 $(x, w)=(0,1)$, and for $ s=4 $, We get the second solution $(x, w)=(119,13)$, with 
  \begin{equation} \label{eq.26} 
   u_{4}+v_{4}\sqrt{2} = 17 + 12 \sqrt{2} ,
\end{equation}
  That is, $ m=12$ and $ n=5$, This is the generating pair of \eqref{eq.05} for the only solution proposed by Fermat problems $ n^{\circ} $ 25, p 340 \cite{2}, using the right triangle whose sides are $ x $ and $ x+1 $.

\end{proof}

\begin{theorem}
There are infinitely many primitive Pythagorean triangles whose the sum of the arms are square.
\end{theorem}

\begin{proof} 
Let $ k $ be a positive integer, Pythagoras showed the sides
 \begin{equation} \label{eq.27}
 x=2k+1, \hspace{1cm} y=\dfrac{(2k+1)^2-1}{2}, \hspace{1cm} z=\dfrac{(2k+1)^2+1}{2}
  \end{equation}
  primitive Pythagorean triangles. The equality  \eqref{eq.27} gives
\begin{equation} \label{eq.28}
 x+y=x+\dfrac{x^2-1}{2} = \dfrac{(x+1)^2-2}{2} = f^2
  \end{equation}
  hence
  \begin{equation*} 
  f^2-2\left( \dfrac{x+1}{2} \right)^2 = -1
  \end{equation*}
by formula \eqref{eq.16}, we obtain 
\begin{equation} \label{eq.29}
  u_{2s-1} + v_{2s-1}\sqrt{2} = f + \left(\dfrac{x+1}{2}\right) \sqrt{2}
  \end{equation}
 for some positive integers $ s $. The formula \eqref{eq.29} gives
\begin{equation} \label{eq.30} 
     f = u_{2s-1} \hspace{2cm} \text{and}  \hspace{2cm}  x = 2 v_{2s-1} -1   
\end{equation}
That is, there are infinitely many primitive Pythagorean triangles whose the 
sum of the arms are square, generated by triples
\begin{equation} \label{eq.31} 
     x = 2 v_{2s-1} -1, \hspace{.3cm}  y = \dfrac{(2 v_{2s-1} -1)^2-1}{2}, \hspace{.3cm}  z = \dfrac{(2 v_{2s-1}-1)^2+1}{2}   
\end{equation} and
\begin{equation*}  
        x+y =  (u_{2s-1})^2     
\end{equation*}
where $ (u_{2s-1}, v_{2s-1}) $ are solutions of \eqref{eq.13} given by the formula \eqref{eq.16}.
\end{proof}

 \begin{table}[h!]
 \begin{tabular}{|c|c|c|c|c|c|c|} 
\hline
$ s $ & $ u_{2s-1} $ & $ v_{2s-1} $ & $ x $ & $ y $ & $ z $ & $ x+y=u_{2s-1}^2 $ \\ \hline
1          &       1     &     1       &     1       &      0      &    1        &       1       \\ \hline
2          &      7      &     5       &      9      &      40      &       41     &         49  \\ \hline
3          &      41      &    29        &     57       &     1624       &   1625         &    1681 \\ \hline
4          &     239       &    169        &    337   &   56784   &     56785       &    57121     \\ \hline
5          &     1393       &    985     &     1969    &    1938480  &  1938481   &       1940449 \\ \hline
 
\end{tabular}
\bigskip
  \caption{The first five primitive Pythagorean triangles whose sum of arms is square, it is given by the formula \eqref{eq.31} (See \href{https://oeis.org/A008843}{A008843}).}
 \end{table}
 
\begin{theorem}
There are infinitely many primitive Pythagorean triangles whose the hypotenuse are square.
\end{theorem}

\begin{proof} 
 Let $ k $ be a positive integer, The equality  \eqref{eq.27} gives
\begin{equation} \label{eq.32}
 z=\dfrac{x^2+1}{2}  = e^2
  \end{equation}
 where $e$ be a positive integer, hence
  \begin{equation*} 
  x^2-2e^2 = -1
  \end{equation*}
by formula \eqref{eq.16}, we obtain 
\begin{equation}\label{eq.33}
  u_{2r-1} + v_{2r-1}\sqrt{2} = x + e \sqrt{2}
  \end{equation}
 for some positive integers $ r $. The formula \eqref{eq.33} gives
\begin{equation*} 
     x = u_{2r-1} \hspace{2cm} \text{and}  \hspace{2cm}  e =  v_{2r-1}     
\end{equation*}
That is, there are infinitely many primitive Pythagorean triangles whose the hypotenuse are square, generated by triples
\begin{equation} \label{eq.34} 
     x =   u_{2r-1},  \hspace{.3cm}  y = \dfrac{( u_{2r-1})^2-1}{2}, \hspace{.3cm}  z = \dfrac{(u_{2r-1})^2+1}{2}   =  (v_{2r-1})^2
\end{equation} 
where $ (u_{2r-1}, v_{2r-1}) $ are solutions of \eqref{eq.13} given by the formula \eqref{eq.16}.
\end{proof}

\begin{table}[h!]
  \begin{tabular}{|c|c|c|c|c|c| } 
\hline
$ r $ & $ u_{2r-1} $ & $ v_{2r-1} $ & $ x $ & $ y $ & $ z = v_{2r-1}^2$   \\ \hline
1          &       1     &     1       &     1       &      0      &    1                \\ \hline
2          &      7      &     5       &      7      &      24      &       25        \\ \hline
3          &      41      &    29        &     41       &     840       &   841          \\ \hline
4          &     239       &    169        &   239   &   28560   &     28561           \\ \hline
5          &     1393       &    985     &     1393    &    970224  &  970225     \\ \hline
 
\end{tabular}
\bigskip
  \caption{The first five primitive Pythagorean triangles whose the hypotenuse are square, given by the formula \eqref{eq.34} (See \href{https://oeis.org/A008844}{A008844}).}
 \end{table}

\section{Primitive Pythagorean Triangles Generators}

Let $ m,n $ are non-zero integers, with $ m>n , m+n \equiv 1 \pmod 2 $ and $ \gcd(m,n)=1 $.
We define the generators of primitive Pythagorean triangles as follows:
\begin{equation*}
 X=m^2-n^2,  \hspace{2cm}  Y=2mn, \hspace{2cm}  Z=m^2+n^2         
\end{equation*}
That is, $Z$ is strictly positive, $Y$ is exactly positive if and only if $m$ and $n$ have the same sign, and $X$ is exactly positive if and only if $ \vert m \vert > \vert n \vert $.

In this part, we are concerned with the negative values of the numbers $ X $ and $ Y $ in addition to the positive values, so that the sum $ X + Y $ is a perfect square. This necessarily cannot be negative at the same time. Likewise, we can infer the values of $ m $ and $ n $ in terms of the values of $ X $ and $  Z$ as follows

\begin{equation} \label{*} 
 m= \; \delta \; \sqrt{\dfrac{Z + X}{2} }    ,  \hspace{.5 cm}    n= \; \lambda \; \sqrt{\dfrac{Z-X}{2} }  , \hspace{.5cm}        
  D= \sqrt{Z}    ,  \hspace{.5cm}    C= \sqrt{X + Y}     
 \end{equation}
Where   
\[
  \delta = \begin{cases}
    -1 & \text{if $X < 0$ and $ Y>0 $} \\
    1 &  \text{otherwise}  
    \end{cases}
  \hspace{.5cm} \text{and}  \hspace{.5cm}
    \lambda = \begin{cases}
     1 & \text{if $X > 0$ and $ Y>0 $} \\
    -1 &  \text{otherwise}  
    \end{cases}
 \]

 For example, if $ m=-5 $ and $ n=-12 $, we get the primitive Pythagorean triples
 \begin{equation*}
 X=-119,  \hspace{1cm}  Y=120, \hspace{1cm}  Z=169,   \hspace{1cm}  C=13, \hspace{1cm}  D=1       
\end{equation*}
 Moreover, if we take
 \begin{equation*}
 X=4565486027761  ,  \hspace{1cm}  Y= 1061652293520   , \hspace{1cm}  Z= 4687298610289        
\end{equation*}
  from \eqref{*} we conclude  
 \begin{equation*}
    m=2150905 ,  \hspace{1cm}  n= 246792,\hspace{1cm}  D=  2165017 ,  \hspace{1cm}    C  =  2372159   
 \end{equation*}

We suppose for a rational number $t$,
 \begin{equation} \label{eq.35}
 x = (t + m)^2 - n^2, \hspace{1cm} y = 2 n(t + m), \hspace{1cm}  z = (t + m)^2 + n^2
 \end{equation}
We have the identity
 \begin{equation*}
\left[(t + m)^2 - n^2 \right]^2 + \left[2 n(t + m)\right]^2 = \left[(t + m)^2 + n^2\right]^2
  \end{equation*}
In order that the numbers \eqref{eq.35} satisfy equation \eqref{eq.04} it is necessary and sufficient that
 \begin{equation*}
(t + m)^2 + n^2 = e^2, \hspace{1cm} (t + m)^2 - n^2 + 2 n(t + m) = f^2
  \end{equation*}
That is,
\begin{equation} \label{eq.36}
t^2 + 2 m t + m^{2}+ n^2 = e^2
 \end{equation}
and
\begin{equation} \label{eq.37}
t^2 + 2(n+m)t + m^2-n^2+2mn = f^2
 \end{equation}
 If  $ m^2-n^2+2mn = C^2 $ and $ m^2+n^2= D^2 $, It becomes
 \begin{equation} \label{eq.38}
t^2 + 2(n+m)t + C^2 = f^2, \hspace{1cm} \text{and} \hspace{1cm}  t^2 + 2 m t + D^{2} = e^2
 \end{equation}
 
\begin{definition}
We call the positive primitive solution of the system of two Diophantine equations \eqref{eq.04} any solution in which all the unknowns are strictly positive, with $ \gcd(x,y)=1 $. For example, Fermat found the smallest solutions  
\begin{equation*}
x_{1}=4565486027761, \hspace{.1cm} y_{1}=1061652293520, \hspace{.1cm} e_{1}=2165017, \hspace{.1cm} f_{1}=2372159
\end{equation*}

We call the negative primitive solution  of the system of two Diophantine equations \eqref{eq.04} any solution in which all the unknowns are strictly positive, except $x$ or $y$, with $ x+y>0 $ and $ \gcd(x,y)=1 $.
For example, Fermat found the smallest solutions  
\begin{equation*}
x=-119, \hspace{1cm} y=120, \hspace{1cm} e=13, \hspace{1cm} f=1
\end{equation*}

\end{definition} 
 
 The importance of addressing the negative primitive solution to system \eqref{eq.04}, that is, the possibility of generating a positive primitive solution to equation \eqref{eq.04}, and there are four possible possibilities for generating solutions to system \eqref{eq.04}, which are as follows:
 \begin{itemize}
 \item  A negative primitive solution generates a positive primitive solution.
\item A negative primitive solution generates a negative primitive solution.
\item A positive primitive solution generates a negative primitive solution.
\item  A positive primitive solution generates a positive primitive solution.
 \end{itemize}
 
 \begin{remark} In general, we obtain the new primitive solution after reduction, that is to say by dividing the solution generated by $ \gcd(x,y) $.
 \end{remark}

 \begin{theorem} 
 Consider the system of two Diophantine equations with the rational unknown $ t $
  \begin{equation} \label{eq.39} 
  t^2+at+c^2= \alpha^{2}          , \hspace{2cm}        t^2+bt+d^2= \beta^{2}
\end{equation}
 where $a,b,c,d,\alpha,\beta$ are relative integers, such that $ cd (d \pm c)(bc \pm ad)\neq 0 $.
  the trivial solution is $ t_{0} = 0 $, moreover the two possible solutions given by:
 \begin{align*}
   t_{1} = \dfrac{ (ad^2-bc^2)^2 -(2(d+c)cd)^2 }{4cd (d+c)(bc+ad)}    \\                
t_{2} = \dfrac{ (ad^2-bc^2)^2 -(2(d-c)cd)^2 }{4cd (d-c)(bc-ad)} 
 \end{align*}
   \end{theorem}
 
 \begin{proof}
  Looking  by the ordinary method of Fermat, we obtain: \\
  The equality \eqref{eq.39} gives
  \begin{equation} \label{eq.40} 
  d^2 t^2+ad^2t+(cd)^2= (\alpha d)^{2}          , \hspace{2cm}        c^2 t^2+bc^2t+(cd)^2= (\beta c)^{2}
\end{equation}
  The equality \eqref{eq.40} gives
  \begin{equation} \label{eq.41} 
   (\alpha d)^{2}  - (\beta c)^{2} =(d^2-c^2) t^2+ (ad^2-bc^2)t 
\end{equation}
  That is 
  \begin{equation} \label{eq.42} 
  (d \alpha-c\beta)( d \alpha + c\beta)  = (d+c)t    \left( (d-c)t+ \dfrac{ad^2-bc^2}{d+c}\right)
\end{equation}
  Or
  \begin{equation} \label{eq.43} 
  (d \alpha-c\beta)( d \alpha + c\beta)  = \left(d-c\right)t   \left( (d+c)t+ \dfrac{ad^2-bc^2}{d-c}\right)
\end{equation}
  The equality \eqref{eq.42} gives
\begin{equation} \label{eq.44} 
   d \alpha-c\beta = (d-c)t+ \dfrac{ad^2-bc^2}{d+c} ,  \hspace{1cm}   d \alpha + c\beta   = (d+c)t    
\end{equation}
by adding and subtracting the two terms of formula  \eqref{eq.44}, we have
\begin{equation} \label{eq.45} 
    \alpha = t + \dfrac{ad^2-bc^2}{2d(d+c)} ,  \hspace{2cm}  \beta   = t - \dfrac{ad^2-bc^2}{2c(d+c)}  
\end{equation}
Substituting into formula \eqref{eq.39}, we obtain the first solution
\begin{equation} \label{eq.46} 
     t_{1} =    \dfrac{ (ad^2-bc^2)^2- \left( 2cd(d+c)\right)^2 }{4cd(d+c)(bc+ad)}  
\end{equation}
The same thing, the equality \eqref{eq.43} gives
\begin{equation} \label{eq.47} 
   d \alpha-c\beta = (d-c)t,  \hspace{1cm}   d \alpha + c\beta   = (d+c)t + \dfrac{ad^2-bc^2}{d-c}    
\end{equation}
by adding and subtracting the two terms of formula  \eqref{eq.47}, we have
\begin{equation} \label{eq.48} 
    \alpha = t + \dfrac{ad^2-bc^2}{2d(d-c)} ,  \hspace{2cm}  \beta   = t + \dfrac{ad^2-bc^2}{2c(d-c)}  
\end{equation}
Substituting into formula \eqref{eq.39}, we obtain the second solution
\begin{equation} \label{eq.49} 
     t_{2} =    \dfrac{ (ad^2-bc^2)^2- \left( 2cd(d-c)\right)^2 }{4cd(d-c)(bc-ad)}  
\end{equation}
For solutions $ t_{1} $ and $ t_{2} $ to be well defined, the denominator must be non-zero, that is:
\begin{equation} \label{eq.50} 
   cd (d \pm c)(bc \pm ad)\neq 0     
 \end{equation}
 \end{proof}
 
 \begin{example} we take:  
 \[ a=34,\hspace{.8cm} b=10, \hspace{.8cm} c=1,\hspace{.8cm} d=13,\hspace{.8cm} m=5,\hspace{.8cm} n=12 \]
 We obtain the first solution $ t_{1}=\frac{2048075}{20566} $ by substituting into expression
  \eqref{eq.35}, generating the positive primitive solution given by Fermat
\begin{equation*}
 x_{1}= 4565486027761 ,  \hspace{1cm}    y_{1}= 1061652293520,   \hspace{1cm}   z_{1}= 4687298610289    
\end{equation*}
\begin{equation*}
 e_{1}= \sqrt{z}  =  2165017 ,  \hspace{2cm}    f_{1} = \sqrt{x + y}  =  2372159   
 \end{equation*}
 
 We obtain the second solution $ t_{2}=-\frac{1582}{13} $ by substituting into expression \eqref{eq.35}, generating the negative primitive solution  
\begin{equation*}
 x= 2276953 ,  \hspace{1.5cm}    y= -473304,   \hspace{1.5cm}   z= 2325625    
\end{equation*}
\begin{equation*}
 e= \sqrt{z}  =  1525 ,  \hspace{2cm}    f = \sqrt{x + y}  =  1343   
 \end{equation*}
 \end{example}
 
 \begin{example} we take: \[ a=2722,\hspace{.3cm} b=3034, \hspace{.3cm} c=1343,\hspace{.3cm} d=1525,\hspace{.3cm} m=1517,\hspace{.3cm} n=-156 \]
 
 We obtain the first solution $ t_{1}=-\frac{6620057924551204}{9320377661925} $ by substituting into expression \eqref{eq.35}, generating the negative primitive solution  
\begin{align*}
  x= 54420629434406206268103685648441 , \hspace{2cm}  \\
      y= -21864804036399372236043874332600,  \hspace{2cm} \\
        z= (7658246457672229)^2  ,   \hspace{2cm}    x+y  =  (5705771236038721)^2 
\end{align*}

 We obtain the second solution $ t_{2}=-\frac{5135764561703}{3240054650} $ by substituting into expression \eqref{eq.35}, generating the negative primitive solution  
\begin{equation*}
 x= -206813120469783031691591 ,  \hspace{1cm}    y= 223005565123008949772400,     
\end{equation*}
\begin{equation*}
 \gcd(x,y)=169,  \hspace{1cm} z= (551491888597)^2,  \hspace{1cm}   x+y =   (127249536947)^2
  \end{equation*}
 \end{example}

 \begin{example} we take: 
 \begin{align*} a=4795394,\hspace{2cm} b=4301810, \hspace{2cm} c=2372159, \\
  d=2165017,\hspace{2cm}  m=2150905 , \hspace{2cm} n= 246792
  \end{align*}
  
 We obtain the first solution \[ t_{1}=\frac{-26417915836091201852959675}{23371832231633450697 }\] by substituting into expression \eqref{eq.35}, generating the positive solution  
\begin{align*}
  x= 535680495873241092745209200472728654679243416759524, \hspace{2cm}  \\
      y=   275163562553703711131840159506694122064526245693280        ,  \hspace{2cm} \\
        z= (24540165251772609623297974)^2        ,   \hspace{4cm} \\
         x+y=(28475323675543089493990702)^2      ,  \hspace{4cm}   \\
          \gcd(x,y)= 2502724 \hspace{6cm}          
\end{align*}
\emph{ After reduction, we conclude the positive primitive solution } 
\textcolor{blue}{ 
\begin{align*}
  x= 214038981475081188634947041892245670988588201, \hspace{2cm} & \\
      y=   109945628264924023237017010068507003594693720 ,  \hspace{2cm} & \\  \hspace{2cm}
             z=   240625698472667313160415295005368384723483849      ,   \hspace{2cm} & \\
        z= (15512114571284835412957)^2                   ,   \hspace{4cm} & \\
         x+y=(17999572487701067948161)^2        , \hspace{4cm}  \\
           \gcd(x,y,z)= 1 \hspace{6cm} &
\end{align*}}

\emph{ Note that this solution has already been found by G. Jacob Martens in a different way \cite{6}.
} 

 We obtain the second solution \[ t_{2}=\frac{-43696202356681630705292570379109}{21491187664558271841998066} \] by substituting into expression \eqref{eq.35}, generating the negative solution  
\begin{align*}
  x=  -21733496857139479252342798266058042991635616764351673608324343 , \hspace{2cm}  \\
      y=   26830078589973982647877216088059498465111507506193928641585824         ,  \hspace{2cm} \\
     z= 34528220382514241522041637315964382769750103796188554357775625 ,  \hspace{2cm} \\ 
        z= (5876071849672554975815554150525)^2        ,   \hspace{4cm} \\
         x+y=( 2257561014199727757943672268941 )^2      ,  \hspace{4cm}  \\
           \gcd(x,y)= 169 \hspace{7cm}          
\end{align*}
we conclude the negative primitive solution 
\begin{align*}
  x=  -128600573119168516286052060745905579832163412806814636735647 , \hspace{2cm}  \\
      y=   158757861479135991999273467976683422870482292936058749358496         ,  \hspace{2cm} \\
     z= 204308996346238115515039274058960844791420732521825765430625 ,  \hspace{2cm} \\ 
        z= (452005526897888844293504165425)^2        ,   \hspace{4cm} \\
         x+y=( 173658539553825212149513251457 )^2      ,  \hspace{4cm}  \\
           \gcd(x,y)= 1 \hspace{7cm}\\         
\end{align*}

 \end{example}

\begin{example} we take: 
 \begin{align*} a=37447100114032399472674,  \hspace{1cm}  b=30155088457762762855210,\\
  c=17999572487701067948161, \hspace{1cm}  d=15512114571284835412957,\\
   m=15077544228881381427605 , \hspace{1cm}   n= 3646005828134818308732
  \end{align*}
  
Let \[  t = \dfrac{ (ad^2-bc^2)^2 -(2(d+c)cd)^2 }{4cd (d+c)(bc+ad)} =\frac{p}{q } \hspace{1cm} \text{with}\hspace{1cm} \gcd(p,q)=1\] 
 by substituting into expression 
 \begin{equation*}  
 X = \left((t + m)^2 - n^2\right)q^2, \hspace{.7cm} Y = \left(2 n(t + m)\right)q^2, \hspace{.7cm}  Z = \left((t + m)^2 + n^2\right)q^2
 \end{equation*}
 If $ \lambda=\gcd(X,Y) $ , then $ (x,y,z)= \left( \dfrac{X}{\lambda} , \dfrac{Y}{\lambda}, \dfrac{Z}{\lambda} \right)  $ is primitive solution of equations
 \begin{equation*}
 x^2+y^2=z^2=u^4, \hspace{2cm}  x+y=v^2
 \end{equation*}
    
  \end{example}

\begin{remark}
Unfortunately, we were unable to complete the calculations, but we hope to obtain a third positive primitive for the problem. 
\end{remark}

\section{Conclusion}

Fermat consider $ m_1=5,n_1=12 $, that is  $ t_1=\dfrac{2048075 }{20566} $, $ q_1=20566^2 $, and conclude the first solution
\[( x_1,  y_1, z_1 )= (4565486027761,1061652293520,4687298610289 ).\]

In general, the algorithm for finding primitive solutions can be generalized recurrently as follows:

  For $ s\geq2 $ , we take 
   \[  m_s=\sqrt{ \dfrac{z_{s-1} +x_{s-1 } }{2}} \hspace{1cm}  \text{and} \hspace{1cm} n_s=\sqrt{ \dfrac{z_{s-1}-x_{ s-1} }{2}}\] and conclude
 \[   t_s=  \dfrac{(ad^2-bc^2 )^2-(2cd(d+c))^2 }{ 4cd(d+c)  (bc+ad)}  = \dfrac{p_{s}}{q_{s}} ,\]
With
\[a=2(m_s+n_s ),\hspace{1cm} b=2m_s  , \hspace{1cm}    c=\sqrt{x_{s-1}+y_{s-1} },\hspace{1cm} d=\sqrt{  z_{s-1}}\]
and $ q_s $ the denominator of irreducible fraction of $ t_s $. By substituting into expression 
 \begin{equation*}  
 X_{s} = \left((t_{s} + m_{s})^2 - n_{s}^2\right)q_{s}^2, \hspace{.3cm} Y_{s} = \left(2 n_{s}(t_{s} + m_{s})\right)q_{s}^2, \hspace{.3cm}  Z_{s} = \left((t_{s} + m_{s})^2 + n_{s}^2\right)q_{s}^2
 \end{equation*}
 If $ \lambda_{s}=\gcd(X_{s},Y_{s}) $ , then $ (x_{s},y_{s},z_{s})= \left( \dfrac{X_{s}}{\lambda_{s}} , \dfrac{Y_{s}}{\lambda_{s}}, \dfrac{Z_{s}}{\lambda_{s}} \right)  $ is primitive solution of equations
 \begin{equation*}
 x^2+y^2=z^2=e^4, \hspace{2cm}  x+y=f^2
 \end{equation*}
    
To obtain the desired solution, it is necessary $  t_s+m_s>n_s$ , otherwise, we say that the solution is negative ($ x_{s} $ or $ y_{s} $ is negative).

\section{Another method to solve Fermat's problem}
  
We are currently reviewing Paul's comment on a question asked by Benvitalis on the website 
\footnote{ \href{https://benvitalenum3ers.wordpress.com/2014/04/03/puzzle-pythagorean-triples-sums-of-sides-squares-part-3/}{Puzzle: Pythagorean triples - sums of sides, Squares - Part 3}
Posted on April 3, 2014 by benvitalis.}, improving it to enrich the article since it falls under the same content.

A little bit of analysis for this one.

If we set $  x = p^2 - q^2$ and $  y = 2pq$ then $  x^2 + y^2$ is equal to the square $  (p^2 + q^2)^2$, and we can make $  p^2 + q^2$ a square by setting $  p = r^2 - s^2$ and $  q = 2rs$, which gives $  x^2 + y^2 = (r^2 + s^2)^4$. This also gives \[  x + y= r^4 + 4r^3 s - 6r^2 s^2 - 4rs^3 + s^4.\]
If we divide through the right hand side of that last polynomial by $  s^4$ and then set $ t_{0} = r_{0}/s_{0}$ we get
$  t_{0}^4 + 4t_{0}^3 - 6t_{0}^2 - 4t_{0} + 1$. This is to be a rational square for some number $ t_{0} $.
Since the leading and trailing coefficients are square we can almost factor it. We find that this polynomial differs from
$  (t_{0}^2 - 2t_{0} + 1)^2$   by $4t_{0}^2 (2t_{0} - 3)$. This gives a rational solution $  t_{0} = 3/2$.
However this would give a negative value of $ x $, so even though it is a solution it isn't the one we require. To do this we need to set $  t_{1} = T_{1} + t_{0}$ in the above polynomial where $ t_{0} $ is an arbitrary rational number. This equates to
\[
T_{1}^4 + (4t_{0}+4) T_{1}^3 + (6t_{0}^2+12t_{0}-6) T_{1}^2 + (4t_{0}^3+12t_{0}^2-12t_{0}-4) T_{1} +  t_{0}^4 + 4t_{0}^3 - 6t_{0}^2 - 4t_{0} + 1
 \]
and we know that $  t_{0} = 3/2$ substituted in reduces to $  T_{1}^4 + 10T_{1}^3 + (51/2) T_{1}^2+ (37/2)T_{1} + 1/16$
This again differs from this square $  (T_{1}^2 - 37T_{1} - 1/4)^2 $  by $T_{1}^2 (84T_{1} - 1343)$. where $  T_{1} = 1343/84$ is a solution, therefore we have $  t_{1} = T_{1} + t_{0}$ which is $  t_{1} = 1343/84 + 3/2 = 1469/84$.
so we can take $  r_{1} = 1469$ and $  s_{1} = 84$ and substitute in
$  p_{1} = r_{1}^2 - s_{1}^2$ and $  q_{1} = 2r_{2}s_{1} $ to give $  p_{1} = 2150905$ and $  q_{2} = 246792$,
Then $  x_{1} = p_{1}^2 - q_{1}^2$ and $  y_{1} = 2p_{1}q_{1}$ giving $  x_{1} = 4565486027761$ and $  y_{1} = 1061652293520$
and finally $  z_{1} = 4687298610289$.

Assuming $ r_1=1469 $, $ s_1=84 $ and defined recurrently for $ k\geq2 $ the following expressions
$  T_{k}^4 + a_{k} T_{k}^3 + b_{k} T_{k}^2 + c_{k} T_{k} + d_{k} $, 
where 
\begin{align*} 
 & a_{k} =  \dfrac{4}{s_{k-1}}(r_{k-1} + s_{k-1}),  \\
  &   b_{k} =  \dfrac{1}{s^{2}_{k-1}}(6 r_{k-1}^{2}+12 r_{k-1}s_{k-1}-6 s^{2}_{k-1}),  \\
  &   c_{k} =  \dfrac{4}{s^{3}_{k-1}}(r^{3}_{k-1}+3 r^{2}_{k-1}s_{k-1}-3 r_{k-1}s^{2}_{k-1}-s^{3}_{k-1}), \\      
  &        d_{k} =  \dfrac{1}{s^{4}_{k-1}}(r^{4}_{k-1}+4 r^{3}_{k-1}s_{k-1}-6r^{2}_{k-1}s^{2}_{k-1}-4r_{k-1}s^{3}_{k-1}+s^{4}_{k-1}).
\end{align*}

That is
\begin{align*} 
 & T_{k}^4 + a_{k} T_{k}^3 + b_{k} T_{k}^2 + c_{k} T_{k} +    d_{k} -\left(T_{k}^2-\alpha_{k} T_{k} -\beta_{k} \right)^2  \\
& =  \left(a_{k} -2\alpha_{k} \right) T_{k}^{3} + \left(b_{k} - \alpha_{k}^{2} + 2\beta_{k} \right) T_{k}^{2} + \left(c_{k} -2\alpha_{k} \beta_{k} \right) T_{k} + \left(d_{k} -\beta_{k}^{2} \right) \\
& = T_{k}^{2} \left[ \left(a_{k} +2\alpha_{k} \right) T_{k} + \left(b_{k} - \alpha_{k}^{2} + 2\beta_{k} \right) \right] \\
& = T_{k}^{2} \left[ \left(a_{k} + \dfrac{c_{k}}{ \sqrt{d_{k}}} \right) T_{k} + \left(b_{k} - \dfrac{c_{k}^{2}}{4 d_{k}}  + 2\sqrt{d_{k}} \right) \right] = 0
\end{align*}
It is by taking 
\[ \beta_{k} =\sqrt{d_{k}},   \hspace{1cm} and \hspace{1cm} \alpha_{k}=\dfrac{c_{k}}{2\beta_{k}}=\dfrac{c_{k}}{2\sqrt{d_{k}}}
\] 
We obtain
\[ T_{k}=     \left(   \dfrac{c_{k}^{2}}{4 d_{k}} - b_{k} - 2\sqrt{d_{k}} \right)/ \left(   a_{k} +\dfrac{c_{k}}{  \sqrt{d_{k}}}  \right)     
\]
Therefore we have $  t_{k} = T_{k} + \dfrac{r_{k-1}}{s_{k-1}}$ which is $  t_{k} = \dfrac{r_{k}}{s_{k}}$, and substitute in
$  p_{k} = r_{k}^2 - s_{k}^2$ and $  q_{k} = 2r_{k}s_{k} $.  
Then $  x_{k} = p_{k}^2 - q_{k}^2$ , $  y_{k} = 2p_{k}q_{k}$  and $  z_{k} = p_{k}^2 + q_{k}^2$.   
\begin{problem}
How can we prove that the irreducible fraction $ t_k = \dfrac{r_{k}}{s_{k}}> 1 $, for all values of $  k$?
\end{problem}
 
\begin{example}
\begin{align*} 
 & a_{2} =  \dfrac{4}{s_{1}}(r_{1} + s_{1})= \dfrac{1553}{21},  \\
  &   b_{2} =  \dfrac{1}{s^{2}_{1}}(6 r_{1}^{2}+12 r_{1}s_{1}-6 s^{2}_{1})= \dfrac{2397697}{1176},  \\
  &   c_{2} =  \dfrac{4}{s^{3}_{1}}(r^{3}_{1}+3 r^{2}_{1}s_{1}-3 r_{1}s^{2}_{1}-s^{3}_{1})= \dfrac{3682162385}{148176}, \\      
  &        d_{2} =  \dfrac{1}{s^{4}_{1}}(r^{4}_{1}+4 r^{3}_{1}s_{1}-6r^{2}_{1}s^{2}_{1}-4r_{1}s^{3}_{1}+s^{4}_{1})= \dfrac{5627138321281}{49787136},\\
&  \beta_{2} =\sqrt{d_{2}}= \dfrac{2372159}{7056},   \hspace{1cm} and \hspace{1cm} \alpha_{2}=\dfrac{c_{2}}{2\sqrt{d_{2}}} = \dfrac{3682162385}{99630678} ,\\
& T_{2}=    \left(   \dfrac{c_{2}^{2}}{4 d_{2}} - b_{2} - 2\sqrt{d_{2}} \right)/ \left(   a_{2} +\dfrac{c_{2}}{  \sqrt{d_{2}}}  \right)  = -\dfrac{5632732605275}{619105033092},\\
& t_{2} = T_{2} + \dfrac{r_{1}}{s_{1}} = \dfrac{123672266091}{14740596026},\quad with \quad r_{2}=123672266091, s_{2}=14740596026 \\
&   p_{2}=r_{2}^2-s_{2}^2 = 15077544228881381427605, q_{2}=2r_{2}s_{2}= 3646005828134818308732, \\
& \textcolor{blue}{   x_{2}= p_{2}^2-q_{2}^2 = 214038981475081188634947041892245670988588201,}\\
&\textcolor{blue}{  y_{2}=2p_{2}q_{2} \quad = 109945628264924023237017010068507003594693720,} \\
& \textcolor{blue}{ z_{2}= p_{2}^2+q_{2}^2  \; = 240625698472667313160415295005368384723483849. } 
\end{align*}
\end{example}

\begin{example}
\begin{align*} 
 & a_{3} = \dfrac{276825724234}{7370298013},  \\
  &   b_{3} =   \dfrac{56170650171048599209011}{108642585600863496338},  \\
  &   c_{3} =    \dfrac{2484091554340576859033091015047365 }{800728232781226638144194176394}, \\      
  &        d_{3} =   \dfrac{323984609740005211871964051960752674583281921}{47212845624163809238534594766102957640976},\\
&  \beta_{3} =\sqrt{d_{3}}= \dfrac{17999572487701067948161}{217285171201726992676 },  \\
& \alpha_{3}=\dfrac{c_{3}}{2\sqrt{d_{3}}} = \dfrac{2484091554340576859033091015047365 }{132662213340952648036309005304093} ,\\
& T_{3}=   -\dfrac{437825148963391521638828389137484882137402791039}{98896159763542844418246069612877176065621596668},\\
& t_{3}   = \dfrac{ 165210121963111493378119359635452016808461 }{41690358767495283586912205131644706811652},\\
&    r_{3}=165210121963111493378119359635452016808461, \\
& s_{3}=41690358767495283586912205131644706811652 \\
&   p_{3}=  2555629838490370378141454677149214109425928622015868323491\\
& 7545002092786780958219417,\\
& q_{3}=  1377533851332754072071006726432564262260558482357302220333\\
& 7181166391614008173975144 , \\
& \textcolor{red}{   x_{3}= p_{3}^2-q_{3}^2 =    463364435981466638655721606795402574490785034897} \\
& \textcolor{red}{4431241700702402006711049469339268033653025345762584764923228 }\\ 
& \textcolor{red}{16713177661024808525677009601280034178493297316385999153 } ,\\
&\textcolor{red}{  y_{3}=2p_{3}q_{3} \quad =   70409332279930883377627285381493031367654903017} \\
&  \textcolor{red}{2246816105298509194049219955867268919807417417483602815295710 }\\
& \textcolor{red}{263978203420907943571429370020707878980026592034112342096  }                 , \\
& \textcolor{red}{ z_{3}= p_{3}^2+q_{3}^2  \; =   84288433829499667821222883781519141823442817924 }\\
& \textcolor{red}{9675666444859533148402499591532256435242317441040296075151069} \\
& \textcolor{red}{741561268252863830134704359471711958336290700625845640625 }                 .  
\\
&  \sqrt{z_{3}}= 290324704132286455037121432138325355009852271302457916 \\
& 25262982715784415755764157625, \\
& \sqrt{x_{3}+y_{3}}= 3416808099353511355218046491734686829295873999139\\
& 8355562578195440360113112814117057
\end{align*}
\end{example}

\end{document}